\documentclass[11pt]{amsart}
\usepackage{amsmath,amssymb,amsthm,xspace}
\usepackage{graphicx}
\usepackage[latin1]{inputenc}
\usepackage[english,francais]{babel}
\newtheorem{thm}{Théorème}[section]

\newtheorem{lem}[thm]{Lemma}
\newtheorem{prop}[thm]{Proposition}
\theoremstyle{definition}
\newtheorem{defn}[thm]{Définition}
\theoremstyle{remark}

\numberwithin{equation}{section}

\begin{document}

\title[IET  non topologiquement faiblement mélangeants]
 {Echanges d'intervalles non topologiquement faiblement mélangeants}
\author{Hadda Hmili}

\address{D\'epartement de Math\'emathiques \\ Facult\'e des Sciences de Bizerte \\ 7021
Jarzouna,
Bizerte, Tunisie}


\begin{abstract} Dans cet article, on prouve un critère d'existence 
de fonctions propres continues non constantes
pour les \'echanges d'intervalles, c'est à dire de non mélange
faible topologique. On construit pour tout entier $m > 3$ des
échanges de $m$ intervalles de rang 2 uniquement ergodiques et non
topologiquement faiblement mélangeant répondant ainsi à une question
de Ferenczi et Zamboni dans \cite{sFlZ06}. On construit aussi pour
tout entier pair $m\geq4$ des échanges de $m$ intervalles poss\'edant des valeurs
propres irrationnelles (avec fonctions  propres associées continues donc non
topologiquement faiblement mélangeant) et poss\'edant aussi des valeurs
propres rationnelles (avec fonctions  propres associées continues par
morceaux) et qui sont soit uniquement ergodiques, soit
non minimaux.

\bigskip

\noindent ABSTRACT. In this paper, we prove a criterion for existence of continuous
non constant eigenfunctions for interval exchange transformations, that is for
non topologically weak mixing. We first construct, for any $m>3$, uniquely
ergodic  interval
exchange transformations of $\mathbb Q$-rank 2 with irrational eigenvalues
associated to continuous eigenfunctions, so that are
 not topologically weak mixing, this answers to a question of Ferenczi and
 Zamboni \cite{sFlZ06}. Then, we construct,  for any even number $m\geq 4$,
 interval exchange transformations of $\mathbb Q$-rank 2 with both  irrational eigenvalues
(associated to continuous eigenfunctions)  and non trivial
 rational eigenvalues (associated to piecewise continuous eigenfunctions), 
moreover these examples can be chosen to be either uniquely
ergodic or non minimal.

\end{abstract}
\maketitle

\section*{Historique}

Les  \'echanges d'intervalles ont été introduits, suivant une idée d'Arnold,
par Katok et Stepin \cite{KS67} qui  ont utilisé les \'echanges de 3
intervalles pour construire des familles de transformations avec spectre
continu simple. Ensuite, Keane \cite{mK75} a donné un critère de
minimalité pour les   \'echanges d'intervalles : la condition idoc (voir la
définition 1.6) et a construit \cite{mK77} des  \'echanges
d'intervalles minimaux et non uniquement ergodiques. Veech \cite{wV82}
a prouvé que presque tout
\'echange d'intervalles est uniquement ergodique et a posé
\cite{wV84} de nombreuses
questions sur les propriétés spectrales réalisables par les \'echanges
d'intervalles. Récemment, Avila et Forni \cite{aA2006}  ont montré
que presque tout échange d'intervalles est faiblement mélangeant, précédement,
Nogueira et Rudolph \cite{aNdR97} avaient  montré que presque tout
échange d'intervalles est topologiquement faiblement mélangeant. Ici, à
l'opposé du cas générique  nous construisons des familles d'échanges
d'intervalles avec  spectre non trivial.

\section{Introduction}

 Soient $r> 0$ un entier et $I^{r} = [0,r[$. Soit $m>0$ un entier. On appelle {\it \'echange de $m$ intervalles} de $I^r$ une bijection de $I^r$ sur $I^r$
 pour laquelle il existe une subdivision $a_{0}= 0 < a_{1} < ... < a_{m} = r$ et des réels $
 \delta_{i}$ tels que $T(x) = x+ \delta_{i}$ pour tout $x\in[a_{i-1}, a_{i}[$.

 Le vecteur  $\delta = (\delta_{1},\delta_{2}, \delta_{3},...,\delta_{m})$
 est dit {\it vecteur de translation} de $T$.

 On note  $D(T) = \{a_{1},a_{2},...,a_{m-1}\}$ l'ensemble des {\it
   discontinuit\'es} de $T$.

On note $I_{i}= [a_{i-1}, \ a_{i}[$ , $i= 1,2,3,..., m$ et $J_{j} = [b_{j-1}, \ b_{j}[$
les intervalles $T(I_{i})$ écrits dans l'ordre. Ainsi, on associe à $T$ une permutation $\pi$
de $\{1,...,m\}$ par $T(I_i) = J_{\pi(i)}$.
\medskip

{\unitlength=0,5 mm
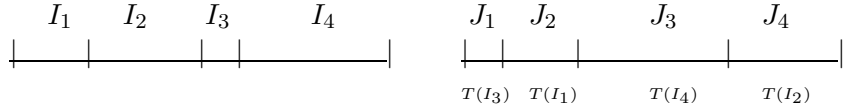
\begin{picture}(250, 25)
\put(0,0){\line(1,0){100}} \put(120,0){\line(1,0){100}}

\put(0,0){$\vert$}
 \put(10,10){$I_1$}
\put(20,0){$\vert$} \put(30,10){$I_2$} \put(50,0){$\vert$}
\put(52,10){$I_3$} \put(60,0){$\vert$} \put(80,10){$I_4$}
\put(100,0){$\vert$}

\put(120,0){$\vert$}
 \put(122,10){$J_1$} \put(120,-10){$\scriptscriptstyle T(I_3)$}
\put(130,0){$\vert$}
\put(138,10){$J_2$}\put(138,-10){$\scriptscriptstyle T(I_1)$}
\put(150,0){$\vert$}
\put(170,10){$J_3$}\put(170,-10){$\scriptscriptstyle T(I_4)$}
\put(190,0){$\vert$}
\put(200,10){$J_4$}\put(200,-10){$\scriptscriptstyle T(I_2)$}
\put(220,0){$\vert$}
\end{picture}}

\bigskip

\bigskip

Dans ce cas on a  $\pi= \left(\begin{array}{ll} 1 \ \ 2 \ \ 3
\ \  4 \\ 2 \ \ 4 \ \ 1 \  \  3\end{array}\right)$
\bigskip

 On associe \`a $T$  son {\it  vecteur longueur} not\'e $\lambda = (\lambda_{1},\lambda_{2},
 \lambda_{3},...,\lambda_{m})$, défini par $\lambda_i = \vert I_i\vert$ la longueur de $I_i$.
\medskip

Le vecteur $\delta$ se déduit de $\lambda$ par :

$$\delta_i = \sum_{j= 1,..., \pi(i)}\lambda_{\pi^{-1}(j)}-\sum_{j=1,...i}\lambda_j.$$
\bigskip

On dit que \ $T$ \ est {\it minimal}  si l'orbite \ $O_{T}(x) =
\{T^{n}(x): n\in \mathbb{Z}\}$ de tout point \ $x\in I^{r}$ \ est
dense dans \ $I^{r}$.

Un échange d'intervalles est dit {\it de rang $2$} si les longueurs
associées \ $\lambda_i$ \ appartiennent à un m\^eme \ $\mathbb
Q$-espace vectoriel de rang $2$. On peut montrer (\cite{mdB88}) que
l'on peut toujours se ramener au cas où les \ $\lambda_i$ \
appartiennent à un $\mathbb Q$-espace vectoriel engendré par \ $1$ \
et un nombre irrationnel $\alpha$. Boshernitzan a montré dans
\cite{mdB88} :
\medskip

{\bf Th\'eor\` eme (Boshernitzan)} \cite{mdB88}.{\it  Tout échange
d'intervalles minimal et de rang \ $2$ \ est uniquement ergodique.}
\medskip

Dans ce m\^eme article, Boshernitzan pose aussi la question de la
stabilité des propriétés spectrales pour les échanges d'intervalles
de rang \ $2$.

\bigskip

\bigskip

\bigskip

\bigskip

\centerline {\bf Notions spectrales.}

\smallskip

\medskip

Soit\ $T$  un \'echange d'intervalles sur un intervalle \ $I$.
On note  $\mu$  la mesure de Lebesgue sur  $I$  et
$L^{2}(I,\mu) = \{ f: I \rightarrow \ \mathbb{C}$ \ telle que 
\ $||f||_{2}=\int|f|^{2}d\mu <\infty\}$.

\begin{defn} On d\'efinit un opérateur unitaire par  :
\begin{center}
$U_{T}: L^{2}(I,\mu)\rightarrow L^{2}(I,\mu)$,\\
$f\mapsto foT$
\end{center}

Cet op\'erateur est appel\'e {\it op\'erateur de composition par \
$T$}.
\medskip

 Soit   \ $a\in \mathbb{C}$, une fonction \ $f\in L^{2}([0,1[,\ \mu)$ \ est appelée {\it fonction
propre de \ $U_{T}$}  associée à la {\it valeur propre } $a$  si  $foT = af, $
\ $\mu$-presque partout.
\end{defn}

\bigskip

\noindent {\bf Remarques.} Par d\'efinition, les fonctions propres
de \ $U_{T}$ \ associ\'ees \`a la valeur propre \ $1$ \ sont les
fonctions constantes sur les orbites de \ $T$.

Comme \ $U_{T}$ \ est unitaire ($||U_{T}(f)||_{2}= ||f||_{2}$), les valeurs propres de \ $U_{T}$ \ sont de module \ $1$ \ et
s'écrivent \ $a = \exp(2i\pi\alpha)$ \ avec \ $\alpha\in [0,1[$.

\smallskip

\begin{defn}
\noindent La valeur propre $a$ est dire {\it rationnelle} [resp.  {\it
  irrationnelle}] si $\alpha$ est  rationnel [resp. 
  irrationnel].\end{defn}

\begin{defn} Soit $T$ un échange d'intervalles.

\noindent - On dit que \ $T$ \ est {\it ergodique} par rapport \`a
la mesure de Lebesgue \ $\mu$ \ si les ensembles \ $T$-invariants
sont de \ $\mu$-mesure pleine ou nulle.

\noindent - \ $T$ \ est dit {\it \ $C^{0}$-ergodique} si les fonctions
propres \emph{continues } de \ $U_{T}$ \ associ\'ees \`a la valeur
propre \ $1$ \ sont les constantes.

\noindent - On dit que \ $T$ \ est {\it uniquement ergodique} si 
 la mesure \ $\mu$ \ est la seule  mesure de probabilit\'e \ $T$-invariante.

\noindent - \ $T$ \ est dit {\it faiblement m\'elangeant} si les fonctions
propres dans \ $L^{2}(I,\mu)$ de \ $U_{T}$ \ sont les constantes.

\medskip

\noindent - \ $T$ \ est {\it topologiquement faiblement m\'elangeant} si les
fonctions propres continues de \ $U_{T}$ \ sont les constantes.
\end{defn}

\bigskip

\noindent {\bf Propri\'et\'es.} 

\smallskip

 i) Si  $T$  est minimal alors
\ $T$ \ est  $C^{0}$-ergodique.

\smallskip

ii) Si  $T$  est ergodique  par rapport \`{a} la mesure de
Lebesgue \ $\mu$ \ alors \ $T$ \ est minimal et donc  \ $T$ \ est \ $C^{0}$-ergodique.

\smallskip

iii) Si  $T$ est uniquement ergodique  alors  $T$  est
 ergodique par rapport \` a $\mu$.

\smallskip

iv) Si \ $T$ \ est faiblement m\'elangeant alors \ $T$ \ est
topologiquement faiblement m\'elangeant. La réciproque est fausse:

\begin{quote} {En effet, soit \ $R_{\alpha}$ \ un \'echange de \ $2$ \ intervalles
avec \ $\alpha \notin \mathbb{Q}$. On sait que les fonctions propres
de \ $U_{R_{\alpha}}$ \ sont de la forme \ $g_{k}(x) = Aexp(2i\pi
kx)$, \ $k\in \mathbb{Z}$, o\` u \ $A$ \ est une constante complexe.
Soit \ $h$\ un \'echange de \ $m$ \ intervalles,  posons \ $T = h
\circ R_{\alpha}\circ h^{-1}$. Alors \ $T$ \ est un \'echange d'au
plus \ $2m$ \ intervalles. Les fonctions \ $f_{k} = g_{k}\circ
h^{-1}$, $k\in \mathbb{Z}$, \ sont les fonctions propres de \
$U_{T}$\ dans \ $L^{2}(I,\mu)$ et \ non constantes. Donc, $T$ n'est
pas faiblement m\'elangeant. Mais, \ $T$ \ est topologiquement
faiblement m\'elangeant car si \ $f$ \ est une fonction propre
continue non constante de \ $U_{T}$ \ alors \ $f\circ h$ \ est une
fonction propre de  \ $R_{\alpha}$, donc \ $f$ \ s'\'ecrit : \ $f =
g_{k}\circ h$, \ pour un \ $k\in \mathbb{Z}^{*}$ \ qui est continue
par morceaux et jamais continue dès que $m>2$.}\end{quote}

\bigskip

\centerline {\bf Questions et R\'esultats.}

\bigskip

 Arnoux  a exhibé dans \cite{pA88} un
exemple non trivial d'échange de $7$ intervalles qui possède une
fonction propre continue non constante. De plus, Nogueira et Rudolph
ont prouvé que les échanges de $3$ intervalles sont topologiquement
faiblement mélangeants (voir \cite{aNdR97}). Dans le paragraphe
questions de \cite{sFlZ06}, S. Ferenczi et L.Q. Zamboni
s'interrogent sur l'existence d'échange de $4$ intervalles possédant
des fonctions propres continues. Ici, nous répondons par
l'affirmative à cette question en montrant les
\bigskip

\begin{thm} Pour tout entier \ $m>3$, pour tout $\beta\in [0,1]$ irrationnel, \ il existe des \'echanges \ $T_{\beta}$
  de \ $m$ \ intervalles uniquement
ergodiques vérifiant la conditon {\it idoc} de Keane (voir la définition 1.6)  et non topologiquement faiblement m\'elangeant avec
fonctions propres continues non constantes associ\'ees \`a la valeur propre \ $exp(2i\pi\beta)$.
\end{thm}
\medskip

\begin{thm} Pour tout entier pair \ $m\geq 4$, pour tout $\alpha\in [0,1]$
  irrationnel,\ il existe des \'echanges \ $T_{\alpha}$  de \ $m$
\ intervalles  non topologiquement faiblement
mélangeants, avec fonctions propres continues non constantes
associ\'ees \`a la valeur propre \ $exp(-2i\pi \alpha)$ \ et  possédant aussi
 $exp(4i\pi \frac{1}{m})$  comme valeur propre rationnelle. De plus, ces
 exemples peuvent \^etre choisis  soit uniquement ergodiques, soit non minimaux.
\end{thm}

\bigskip

\noindent {\bf Remarques.} Ce dernier résultat  constitue un phénomène spectral nouveau
par rapport aux rotations. D'autre part, tous les exemples construits pour ces
deux théorèmes sont de rang 2 ; 
dans la proposition 2.1 B, nous montrons que si l'on cherche  des échanges
d'intervalles avec fonctions propres dérivables à dérivée dans $L^2 (I,\mu)$ alors
nécessairement ces  échanges d'intervalles  sont de rang 2.  Par ailleurs, nos
constructions restent valables  pour les valeurs rationnelles des  paramètres
 et fournissent des valeurs propres rationnelles pour des échanges d'intervalles qui ne sont plus minimaux. 

\bigskip

\bigskip

\bigskip

\bigskip

\bigskip

\centerline{\bf Critère de minimalité - Condition idoc de Keane.}

Nous indiquons ici un critère de minimalité dû à Keane \cite{mK75}.

\begin{defn} Une permutation \ $\pi$ \ de \ $\{1, 2,..., m\}$ \ est \emph{irréductible} si

$\pi(\{1, 2,.., k\} ) \neq \{1,2,.., k\}$ pour tout \ $1\leq k < m$.
\medskip

On dit que \ $T$ \ satisfait \emph{la condition idoc} si :
\medskip

a) \ $O(a_{i}) =\{T^{n}(a_{i}): n\in \mathbb{Z}\}$ \ est infini
pour tout $1\leq i< m$,

b) $O(a_{i})\cap O(a_{j}) = \emptyset$\ pour tous \ $1\leq i,
j\leq m-1, \ i\neq j$.
\end{defn}
\medskip

\noindent {\bf Th\'eor\` eme (Keane \cite{mK75}).}  {\it  Si \ $T$ \ satisfait la condition idoc et \ $\pi$ \ est irréductible alors \ $T$ \ est minimal.}

\section{Existence de fonctions propres continues}

\begin{prop} Soit \ $T$ \ un échange de \ $m$ \ intervalles de \ $I$, \ de vecteur longueur
associé \ $(\delta_i)_{1\leq i\leq m}$.
\medskip

A) S'il existe deux réels \ $r$ \ et \ $s$, \ et des entiers \ $p_i$
\ tels que \ $\delta_{i}= r+ p_{i}s$, \ pour tout \ $i\in
\{1,...,m\}$ \ alors \ $T$ \ admet \ $f(x) = \exp(\frac{2i\pi}{s}x)$
\ comme fonction propre associ\'ee \`{a} la valeur propre \
$\exp(2i\pi\frac{r}{s})$.

\smallskip

B) Si \ $T$ \ est uniquement ergodique et admet une fonction propre
dérivable \ $f$ \ dont la dérivée est dans \ $L^2(I,\mu)$ \ alors \
$f(x) = C \ exp(2i\pi k x)$, \ pour un \ $k\in \mathbb{Z}$ \ et \
$T$ \ est de rang \ $2$.
\end{prop}

\noindent {\bf Preuve.}
\smallskip

A) On calcule \ $foT$.
\smallskip

On a \ $f\circ T(x)= \ \exp(\frac{2i\pi}{s}T(x))$, \ pour tout \ $x\in I_{i}$.\\
  \hspace*{3cm}= \  $\exp(\frac{2i\pi}{s}(x+\delta_{i}))$\\
 \hspace*{3cm}= \ $\exp(\frac{2i\pi}{s}(x+r+p_{i}s))$\\
 \hspace*{3cm}= \ $\exp(\frac{2i\pi}{s}x)\exp(2i\pi\frac{r}{s})\exp(2i\pi p_{i})$\\
 \hspace*{3cm}= \ $\exp(2i\pi\frac{r}{s}) f(x)$
\bigskip

B) Par hypothèse, on a  \ $f\circ T = \lambda f$. D\'erivons, on
obtient : \ $f'\circ T = \lambda f'$. Puisque \ $f$ \ n'est pas
constante, il existe \ $x\in I$ \ tel que \ $f(x)\not=0$. Par
cons\'equent \ $f(T^n(x)) = \lambda^n f(x) \not = 0$.

Comme \ $T$ \ est minimal, l'orbite de tout point \ $x$ \ par \ $T$\
est dense. Donc \ $f(y)\not=0$ \ pour tout \ $y \in I$ \ (par
continuit\'e de \ $f$). On a alors:  \ $\frac {f'}{f} \circ T=
\frac{f'}{f}$,\ $\mu$-presque partout.

Par unique ergodicité de \ $T$, \ on obtient :  \ $\frac{f'}{f}= C$,
o\`{u} \ $C$ est une constante complexe non nulle. En intégrant on
obtient \ $f(x) = A exp(C x)$, pour tout $x\in I$, où $A$ est une
constante complexe non nulle.

Ainsi, $\ exp(C T (x) ) = \lambda \ exp(C x)$, \ pour tout \ $x\in
I$. Ecrivons \ $\lambda =\ exp(i \beta)$ \ et \ $C = a+ib$. Par
conséquent,  $\ exp(aT(x)) = \ exp(ax)$ \ d'o\`{u} $a=0$ et $\
exp(ibT(x))= \ exp(i \beta)\ exp(ibx)$ d'o\`{u} $bT(x)= \beta+bx+k(x)$,
avec $k(x)\in \mathbb{Z}$  et finalement  $T (x) = x + \frac{\beta}{b}
+\frac{ k(x)} {b}$. Les $\delta_i$ sont dans le \ $\mathbb Q$-espace
vectoriel engendré par \ $\frac{\beta}{b}$\ et \ $\frac{1} {b}$, \
donc \ $T$ \ est de rang \ $2$.
\bigskip

\section{Preuve du Théorèmes 1.4}

{\bf Construction des  $T_{\beta}$ : échanges de   $m$ intervalles uniquement
 ergodiques vérifiant la conditon {\it idoc} de Keane  et non topologiquement faiblement m\'elangeant.}

\medskip

Soit \ $m>3$ \ un entier et \ $\alpha\in ]0, \ \frac{1}{m-1}[$ \
un irrationnel. Soit \ $T $ \ l'\'echange
d'intervalles de \ $[0,1[$  de permutation associ\'{e}e :  $$\pi = \left (
\begin{array}{ll} \ 1 \ \ \ \ \  \ 2  \ \ \ \ .\ .\ .\ .\ m-2 \ \ \ \ m-1 \ \ \ \ m
 \\ m-1 \ \ \ \ m-2
  \ \ . \ .\ .\ .\ \ 2 \ \ \ \ \ \ m \ \ \ \  1 \end{array}\right)$$
\noindent  et de vecteur longueur associ\'{e} : $$\lambda =
(\frac{1}{m-1},\frac{1}{m-1}, ..,\frac{1}{m-1},
\frac{1}{m-1}-\alpha,\alpha).$$
\bigskip
\medskip

{\unitlength=0,26mm \hskip 2.3 truecm \begin{picture}(300,300)

\put(0,0){\line(1,0){300}} \put(0,300){\line(1,0){300}}
\put(0,0){\line(0,1){300}} \put(300,0){\line(0,1){300}}

\put(40,0){\dashbox(0,300){}} \ \put(80,0){\dashbox(0,300){}}
\put(250,0){\dashbox(0,300){}} \ \put(270,0){\dashbox(0,300)}


 \put(0,280){\dashbox(300,0){}}
\put(0,240){\dashbox(300,0){}}
\put(0,200){\dashbox(300,0){}}
 \put(0,30){\dashbox(300,0){}}

\put(-10,-10){$0$}\put(305,-10){$\scriptstyle 1$}
 \put(30,-15){${\frac{1}{m-1}}$}
  \put(70,-15){${\frac{2}{m-1}}$}

 \put(238,-15){${\scriptscriptstyle \frac{m-2}{m-1}}$}
 \put(269,-15){${\scriptstyle 1-\alpha}$}

 \put(-15,30){$\alpha$}

 \put(-60,200){${\frac{m-4}{m-1}+\alpha}$}
 \put(-60,240){${\frac{m-3}{m-1}+\alpha}$}
 \put(-60,280){${ \frac{  m-2}{ m-1}+\alpha}$}
 \put(-30,300){$ 1$}

\put(0,240){\line(1,1){40}}

\put(40,200){\line(1,1){40}}

\put(250,280){\line(1,1){20}}

\put(270,0){\line(1,1){30}}
\end{picture}}\\\\

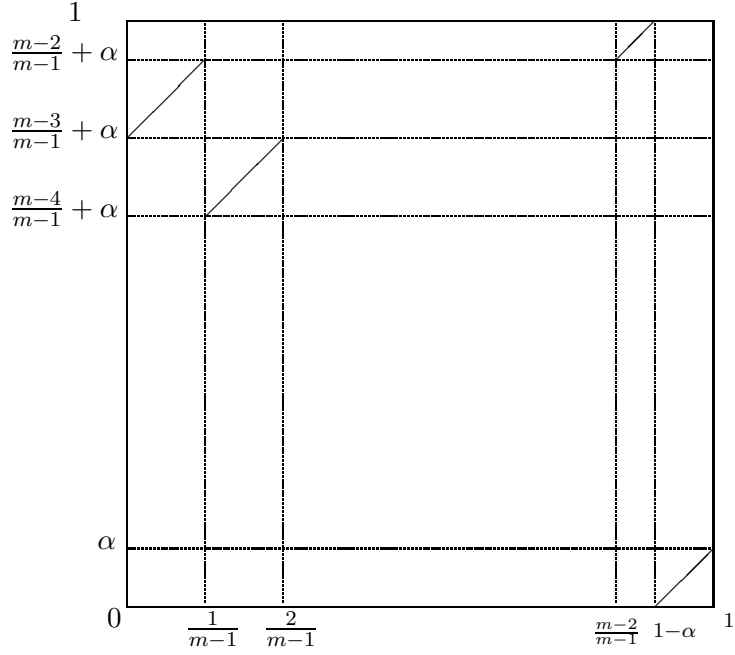
\centerline {\it Figure 1 }

\bigskip
\medskip

On d\'etermine facilement que :
\medskip

- Les discontinuit\'es de \ $T$ \ sont : \ $\frac{1}{m-1},
\frac{2}{m-1}, \ ...\ , \frac{m-2}{m-1}, 1-\alpha$.

\smallskip

- Les images des discontinuit\'{e}s de \ $T$ \ sont :
$\frac{m-4}{m-1}+\alpha, \frac{m-5}{m-1}+\alpha, ....,
\frac{1}{m-1}+\alpha, \alpha, \frac{m-2}{m-1}+\alpha, 0$.

\smallskip

- Les translations de \ $T$ \ sont :

 - $\delta_{i} = \frac{m-1-2i}{m-1}+ \alpha$, pour \ $1\leq i\leq m-2 $

- $\delta_{m-1} = \alpha-\frac{3-m}{m-1}$

- $\delta_{m} = \alpha-1$.

\medskip

Pour tout \ $l\in \mathbb{N}$, \ on a \ $T^{l}(x) =
x+\sum_{i=1}^{m}N_{i}(x,l) \delta_{i}$ \ o\`u \ $N_{i}(x,l) = \sharp
\{0 \leq k \leq l-1: T^{k}x\in I_{i}\}$. Par cons\'equent \
$T^{l}(x) = x+ \frac{p_l(x)}{m-1}+l\alpha$ \ avec \ $p_{l}(x) \in
\mathbb{Z}$.
\medskip

L'échange \ $T$ \ v\'erifie la condition idoc, en effet  :

 - Les orbites des points de discontinuit\'e de \ $T$ \ sont infinies : 

\begin{quote} S'il existe \ $a \in D(T)$ \ et \ $n\in \mathbb{N^{*}}$
 telsque \ $T^{n}(a) = a$, \ alors \ $p_{l}(a)+ n\alpha = 0$, ceci est  absurde par irrationnalité de \ $\alpha$.\end{quote}

 - Les orbites des points de \ $D(T)$ \ sont distinctes: 
\begin{quote} Si deux discontinuités $a_i$ et $a_j$ de \ $T$ \ sont sur une m\^eme
orbite,  alors il existe un entier  $l > 0$ (quitte à échanger $a_i$ et $a_j$)
 tel que :  $T^{l}(a_{i}) = a_{j}$.

\noindent  Alors on aura \ $T^{l}(\frac{i}{m-1}) = \frac{j}{m-1}$  \ ou
 \ $T^{l}(1-\alpha) = \frac{j}{m-1}$ \ ou \ $T^{l}(\frac{i}{m-1}) = 1-\alpha$,
 pour $1\leq i\neq j\leq m-2.$

  Si \ $T^{l}(\frac{i}{m-1}) = \frac{j}{m-1}$ \ alors \ $\frac{p_{l}(\frac{i}{m-1})}{m-1}+l\alpha = \frac{j-i}{m-1}$ \ ce
qui est impossible par irrationnalit\'e de \ $\alpha$.

 Si \ $T^{l}(\frac{i}{m-1}) = 1-\alpha$, \ $1\leq i\neq j\leq
 m-2$ \ alors  \ $\frac{p_{l}(\frac{i}{m-1})}{m-1}+(l+1)\alpha =
1- \frac{i}{m-1}$, \'egalement impossible par irrationnalit\'e de
\ $\alpha$.

 Le cas \ $T^{l}(1-\alpha) = \frac{j}{m-1}$ \ implique
que \ $\frac{p_{l}(1-\alpha)}{m-1} +(l-1)\alpha =
\frac{j}{m-1}-1$, \ n'est possible que si \ $l=1$ \ auquel cas, on
a \ $T(1-\alpha) = \frac{j}{m-1} = 0$ \ ce qui est impossible ($0$ n'est pas  une discontinuité de $T$).\end{quote}

\medskip

Puisque la permutation \ $\pi$ \ est irr\'eductible et que \ $T$ \
v\'erifie la condition idoc alors d'après le théorème de Keane, \ $T$ \
est minimal. Comme les \ $\lambda_{i}$ \ sont dans le \ $\mathbb{Q}$-espace vectoriel engendré par \ $1$ \ et \ $\alpha$ \ alors  \ $T$ \
est de rang \ $2$ \ et donc d'après le th\'eor\` eme de
Boshernitzan, \ $T$ \ est uniquement ergodique.

D'autre part, les \ $\delta_{i}$ \ s'écrivent \ $\delta_{i} = r+
p_{i}s$ avec  $\ r = \alpha$ \ et \ $s = \frac{1}{m-1}$, \ $p_{i}\in
\mathbb{Z}$. Donc d'apr\`es la proposition 2.1 A), \ $T$\ admet
\ $f(x) = \exp(2i \pi (m-1)x)$\ comme fonction propre associ\'ee \`a
la valeur propre \ $\exp(2i\pi (m-1)\alpha)$. On obtient le théorème 1.4 en
posant $\beta = (m-1) \alpha$ et $T_\beta=T$.
\\

\bigskip

{\bf Remarques.} 1) Dans cette famille d'échanges d'intervalles\
$T_{\beta}$ \ de rang \ $2$ \ à un paramètre \ $\beta$, \ on voit
que le non mélange faible topologique est une propriété stable, comme
conjecturé par Boshernitzan dans  \cite{mdB88}. 

2) Pour \ $m = 3$, \ $T$ est un \'echange de \ $2$ \
intervalles et  \ $f(x)=\exp(4i\pi x)$ \ est une fonction
propre associ\'{e}e \`a la valeur propre \ $\exp(4i\pi \alpha)$.
\\

\bigskip

\bigskip

\bigskip

{\bf Cas particuliers.}

$\bullet$ \ \noindent {\bf m =4.} Dans ce cas, \ $T$ \ est un \'echange de \ $4$ \ intervalles de \
$[0, 1[$ 

de permutation associée: \ $\pi = \left (\begin{array}{ll}1 \ \ 2
\ \ 3 \  \
 4 \\ 3 \ \ 2 \ \ 4 \ \ 1 \\
\end{array}\right)$ \ et 

de vecteur longueur: \ $\lambda =
(\frac{1}{3},\frac{1}{3},\frac{1}{3}-\alpha,\alpha)$, \ où \
$\alpha\in ]0,\frac{1}{3}[$ \ est un irrationnel.

\bigskip

\hskip 1.5 truecm

\unitlength=0,3mm

\begin{picture}(180,180)

\put(0,0){\line(1,0){180}} \put(0,180){\line(1,0){180}}
\put(0,0){\line(0,1){180}} \put(180,0){\line(0,1){180}}

\put(-15,-15){$0$}
 \put(60,-20){${\frac{1}{3}}$}
\put(120,-20){${\frac{2}{3}}$}
 \put(140,-20){${1- \alpha}$}
 \put(180,-20){$1$}


 \put(-35,160){${\frac{2}{3}+\alpha}$}
 \put(-35,100){${\frac{1}{3}+\alpha}$}
\put(-15,40){$\alpha$}

 \put(-15,180){$1$}

\put(60,0){\dashbox(0,180){}}
 \put(120,0){\dashbox(0,180){}}
\put(140,0){\dashbox(0,180){}}

\put(0,40){\dashbox(180,0){}} \put(0,100){\dashbox(180,0){}}
\put(0,160){\dashbox(180,0){}}

\put(0,100){\line(1,1){60}}

\put(60,40){\line(1,1){60}}

\put(120,160){\line(1,1){20}}

\put(140,0){\line(1,1){40}}

\end{picture}

\smallskip

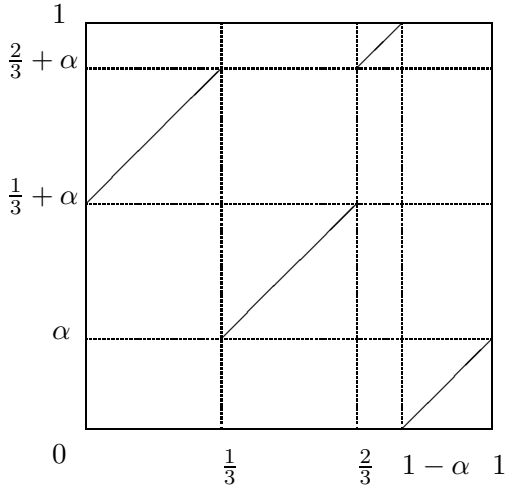
\centerline {\it Figure 2 }

\medskip

On v\'erifie que \ $f(x) = \exp(6i\pi x)$ \ est une fonction propre continue 
associ\'ee \`a la valeur propre \ $\exp(6i\pi \alpha)$. \medskip

$\bullet$ \ \noindent {\bf m = 5.} Dans ce cas, \ $T$ \ est un \'echange de \
$5$ \ intervalles de \ $[0,1[$

 de permutation associ\'{e}e :  $\pi = \left (
 \begin{array}{ll} 1 \ \ 2 \ \ 3 \  \  4  \ \ 5 \\ 4 \ \ 3 \ \ 2 \  \
 5 \ \ 1  \end{array}\right)$ et 

de vecteur longueur associ\'e : \ $\lambda =
(\frac{1}{4},\frac{1}{4},\frac{1}{4},\frac{1}{4}-\alpha,\alpha)$\
où \ $\alpha \in ]0,\ \frac{1}{4}[$ \ est un irrationnel.

Le vecteur de translation \ $\delta$ \ est donn\'e par : \ $\delta =
(\alpha+\frac{1}{2}, \alpha, \alpha - \frac{1}{2}, \alpha,
\alpha-1)$.
\bigskip

\hskip 3 truecm \begin{picture}(200,200)

\put(0,0){\line(1,0){200}} \put(0,200){\line(1,0){200}}
\put(0,0){\line(0,1){200}} \put(200,0){\line(0,1){200}}

\put(-15,-15){$0$}

 \put(50,-20){${\frac{1}{4}}$}

 \put(100,-20){${\frac{1}{2}}$}

 \put(150,-20){$\frac{3}{4}$}

 \put(160,-15){$1-\alpha$}

  \put(200,-20){$1$}

 \put(-15,40){$\alpha$}

 \put(-35,90){$\alpha+\frac{1}{4}$}

 \put(-35,140){$\alpha+\frac{1}{2}$}

 \put(-35,190){$\alpha+\frac{3}{4}$}

 \put(-20,200){$1$}

\put(50,0){\dashbox(0,200){}} \put(100,0){\dashbox(0,200){}}
\put(150,0){\dashbox(0,200){}}\put(160,0){\dashbox(0,200){}}

\put(0,40){\dashbox(200,0){}} \put(0,90){\dashbox(200,0){}}
\put(0,140){\dashbox(200,0){}} \put(0,190){\dashbox(200,0){}}

\put(0,140){\line(1,1){50}}

\put(50,90){\line(1,1){50}}

\put(100,40){\line(1,1){50}}

\put(150,190){\line(1,1){12}}

\put(160,0){\line(1,1){40}}
\end{picture}\\\\

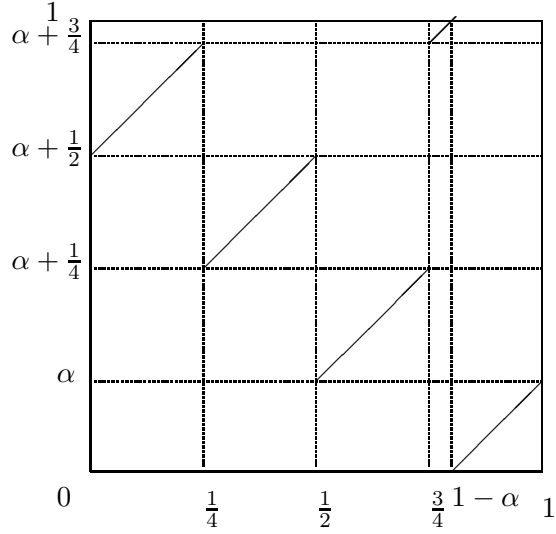
\centerline {\it Figure 3 }

\medskip

On vérifie alors que \ $f(x) = \exp(8i\pi x)$ \ est une fonction
propre de \ $T$ \ associ\'ee \`a la valeur propre \ $\exp(8i\pi
\alpha)$.

\section{Preuve du Th\'eor\` eme 1.5}

{\bf Construction des  $T_{\alpha}$ : échanges de   $m$ 
intervalles  avec valeurs propres irrationnelles et rationnelles 
non triviales  uniquement ergodiques ou non minimaux .}

\medskip

Soient \ $m = 2n$ \ un entier pair, \ $\sigma$\ une permutation de
$\{1,..., n\}$ \ et \ $\alpha\in [0,1]$ \ un irrationnel. Soit \ $T$ \
l'échange de \ $m$ \ intervalles \ de \ $[0,n[$ \ défini  par :

\medskip

 $T([i-1,i[) = [\sigma (i) -1, \sigma (i)[$ pour $i= 1,.., n$ et 
sur $[i-1,i[$ on pose 
$$T(x) = \left\{\begin{array}{lllll}
 x+1-\alpha +\sigma (i)-i & \textrm{si $x\in [i-1,
i-1+\alpha[$}\\
 x-\alpha+\sigma (i)-i & \textrm{si $x\in [
i-1+\alpha,i[$}. 
\end{array}\right.$$
\medskip

\medskip

\noindent Le  vecteur longueur associé  à \ $T$ \ est  $\lambda =
(\alpha,1-\alpha, \alpha,1-\alpha, ....,\alpha,1-\alpha )$. Donc, \
$T$ \ est de rang \ $2$. 

\noindent Les discontinuités de \ $T$ \ sont : \
$\alpha, 1, 1+\alpha, 2,..., (n-1)+\alpha$.

\noindent  Le vecteur translation \
$\delta$ \ est :
$$\delta_{j} = \left\{\begin{array}{lllll}
 1-\alpha +\sigma (i)-i & \textrm{si \ $j=2i-1$ \ est impair}\\
 -\alpha+\sigma (i)-i & \textrm{si \ $j=2i$ \ est pair}; 
\end{array}\right.$$

\noindent où $j=1,...,2n$ et $\ i = 1,..., n$.

\bigskip

La preuve du th\'eor\` eme 1.5 \ est cons\'equence des lemmes 4.1 et
4.2 ci-dessous :
\medskip

\begin{lem} Les trois propriétés suivantes sont équivalentes :

i) \ $T$ \ est uniquement ergodique

ii) \ $T$ \ est minimal

iii) \ $\sigma$ \ est un \ $n$-cycle
\end{lem}
\bigskip

\noindent {\bf Preuve.} $i) \Longrightarrow ii)$ : est clair (voir les
propriétés des notions spectrales).

$ii) \Longrightarrow i)$ : r\'esulte du fait que \ $T$ \ est de
rang \ $2$ \ et du th\'eor\` eme de Boshernitzan.

\medskip

$ii)\Longrightarrow iii)$ : Si \ $\sigma$ \ n'est pas un \ $n$-
cycle,  $\sigma$ contient un cycle de longueur  $p < n$, donc  \ il existe \
$p < n$ tel que \ $T^p ([0,1[) = [0,1[$. Comme : 

\smallskip

$T([0,1[) \ \ \ \ \ = \ [\sigma(1)-1,\sigma(1)[$, 

 $T^2([0,1[)  \ \ \ \  = \ 
 [\sigma^2(1)-1,\sigma^2(1)[$    

 $\vdots$

$T^{p-1}([0,1[) \  = \ [\sigma^{p-1}(1)-1,\sigma^{p-1}(1)[$,  \ \ alors

$\displaystyle \bigcup_{k=0}^{p-1} T^k([0,1[)$ \ est un ensemble invariant de
mesure $p$ donc  distinct de \ $[0,n[$. L'\'echange \ $T$ \ n'est donc pas
minimal.

\medskip

$iii)\Longrightarrow ii)$: \ si \ $\sigma$ \ est un \ $n$-cycle,
alors \ $\displaystyle\bigcup_{k=0}^{n-1} T^k([0,1[) = [0,n[$. L'application de
premier retour de \ $T$ \ sur \ $[0,1[$ \ est la rotation
irrationnelle d'angle \ $n(1-\alpha)$. \ Par suite, \ $T$ \ est
minimal.
\medskip

\begin{lem}
a) \ $f(x) = exp(2i\pi x)$ \ est une fonction propre de \ $U_T$ \ de
valeur propre \ $exp(-2i\pi \alpha)$.

b) \ $exp(2i\pi\frac{1}{ n})$ \ est valeur
propre rationnelle de \ $U_T$.
 \end{lem}

\noindent {\bf Preuve.}

a) On vérifie directement avec la proposition 2.1 A  que $f(x) = \exp(2i\pi x)$ est une fonction
propre associ\'{e}e \`{a} la valeur propre $\exp(-2i \pi \alpha)$.

\medskip

b) On voit que $i$  vérifie $\sigma ^n(i)=i$ donc  on a  $T^n(I_i) = I_i$, où
$I_i=[i-1,i[$  et donc la fonction 

\centerline{ $\displaystyle \phi(x) = exp (2i\pi \left ( \sum_{k=0}^{n-1}
    \frac{k}{n}   \    \mathbb I_{T^k(I_1)} \right) ) $}

\noindent   est une
fonction propre associée à la valeur propre rationnelle $\exp (\frac{2i\pi}
{n})$, où $ \mathbb I_J$ représente la fonction indicatrice de l'ensemble $J$.

 En effet, soit $x\in I_i$, pour $i\in \{1,...,n\}$ par
  définition de $T$ et puisque $\sigma$ est un cycle, il existe $k_i$ tel que
  $I_i= T ^{k_i} (I_1)$, ainsi $x\in {T^{k_i}(I_1)}$ et $T(x) \in {T^{k_i
      +1}(I_1)}$, on a alors : 

\medskip

\centerline{ $\displaystyle \phi(x) = \exp (2i\pi (  \frac{k_i}{n} \ 
   \mathbb I_{T^{k_i}(I_1)} (x)  ) =  \exp (2i\pi \frac{k_i}{n})$ \ \ \ et }

\medskip

\centerline{ $\phi(T(x)) = \exp \bigl( \ \ 2i\pi \ \bigl( \frac{k_i +1}{n}
  \  \mathbb I_{T^{k_i+1}(I_1)}(T(x)) \ \bigl) \ \  \bigl ) \ =  \  \exp (2i\pi \frac{k_i+1}{n}) $}  

\smallskip

\centerline{ $ =  exp (2i\pi
  \frac{1}{n})\exp (2i\pi \frac{k_i}{n}) = \exp (2i\pi
  \frac{1}{n}) \phi(x)$.}

\bigskip

{\bf Exemple 1.}  Soit \ $m = 2n$ \ un entier pair et \ $\alpha$ \ un
irrationnel. Soit  $\sigma$ la permutation donnée par  :  \ $$ \sigma  = \left (
\begin{array}{ll} \ 1 \ \ \ \ \  \ 2  \ \ \ \ .\ .\ .\ .\ \ \ \ n-1 \ \ \ \ n \\\
2\ \ \ \ \ \ 3 \ \ \ . \ .\ .\ . \ \ \ \ \ \ \ \ \  n \ \ \ \ \ 1
  \end{array}\right).$$

On considère  $T$ \ l'\'echange de \ $m$ \ intervalles de \
$[0,n[$ construit comme au début de cette section 4. Puisque \ $\sigma$ \ est un \ $n$-cycle, alors d'apr\` es le lemme 4.1, l'\'echange \ $T$ \ est uniquement ergodique.
 La fonction \ $f(x) = \exp(2i\pi x)$ \ est une fonction propre de \
$U_{T}$ \ associée à la valeur propre \ $\exp(-2i\pi \alpha)$, d'après le
lemme 4.2 a. Pour \ $m = 12$, on a la figure ci-dessous:\\
\medskip

\bigskip

{\unitlength=0,23mm \hskip 3 truecm \begin{picture}(300,300)
\put(0,0){\line(1,0){300}} \put(0,300){\line(1,0){300}}
\put(0,0){\line(0,1){300}} \put(300,0){\line(0,1){300}}

\put(50,0){\line(0,1){300}}
 \put(100,0){\line(0,1){300}}
 \put(150,0){\line(0,1){300}}
 \put(200,0){\line(0,1){300}}
\put(250,0){\line(0,1){300}}

 \put(0,50){\line(1,0){300}}
 \put(0,100){\line(1,0){300}}
 \put(0,150){\line(1,0){300}}
 \put(0,200){\line(1,0){300}}
 \put(0,250){\line(1,0){300}}

\put(-10,-10){$0$}  \put(19,0){$\scriptstyle \vert$}
\put(18,-12){$\scriptstyle \alpha$}
\put(47,-12){$1$} 
 \put(69,0){$\scriptstyle \vert$} \put(63,-16){$\scriptstyle 1+\alpha$}\put(97,-12){$2$}
\put(147,-12){$3$}
 \put(197,-12){$4$} \put(247,-12){$5$} \put(300,-12){$6$}

 \put(-15,50){$1$}  \put(-15,100){$2$}\put(-15,150){$3$}\put(-15,200){$4$}
 \put(-15,250){$5$} \put(-15,300){$6$}

\put(0,80){\line(1,1){20}}\put(20,50){\line(1,1){30}}
\put(50,130){\line(1,1){20}}\put(70,100){\line(1,1){30}}
\put(100,180){\line(1,1){20}}\put(120,150){\line(1,1){30}}

\put(20,50){\dashbox(0,50){}} \ \put(0,80){\dashbox(50,0){}} \
\put(70,100){\dashbox(0,50){}} \ \put(50,130){\dashbox(50,0){}}
\put(120,150){\dashbox(0,50){}} \put(100,180){\dashbox(50,0){}} \

\put(150,230){\line(1,1){20}}\put(170,200){\line(1,1){30}}
\put(200,280){\line(1,1){20}}\put(220,250){\line(1,1){30}}
\put(250,30){\line(1,1){20}}\put(270,0){\line(1,1){30}}

\put(170,200){\dashbox(0,50){}} \ \put(150,230){\dashbox(50,0){}}
\ \put(220,250){\dashbox(0,50){}} \
\put(200,280){\dashbox(50,0){}} \put(270,0){\dashbox(0,50){}} \
\put(250,30){\dashbox(50,0){}}\end{picture}} \\\\
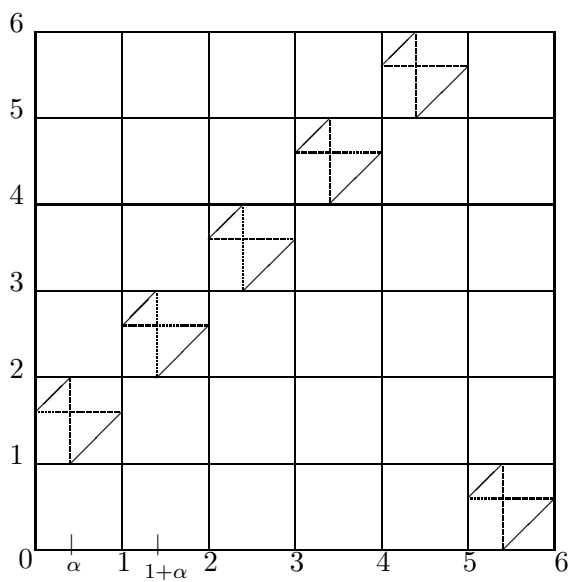
\centerline{\it Figure 4}

\bigskip

{\bf Exemple 2}

\medskip

 \noindent Soit \ $m = 2n$ ($n\geq 1$) \ et \ $T$ \ l'\'echange de \ $m$ \
intervalles de \ $I^{n}=[0,n[$ correspondant  \` a la permutation  \ $ \sigma
= \left (\begin{array}{ll} \ 1 \ \ \ \ \  \ 2  \ \ \ \ .\ .\ .\ .\ \ \ \ n-1 \ \ \ \ n \\\
n \ \ \ \ \ \ n-1 \ \ \ \ .\ .\ .\ .\ \ \ \  2 \ \ \ \ \ 1
  \end{array}\right)$.

\bigskip

- La permutation \ $\pi$ \ associ\'ee \` a  \ $T$ \ est : $\pi = \left (
\begin{array}{ll} \ 1 \ \ \ \ \  \ 2  \ \ \ \ .\ .\ .\ .\ \ \ m-1 \ \ \ \ m
 \\ m \ \ \ \ m-1
  \ \ . \ .\ .\ .\ \ 2 \ \ \ \ \ \ \ \ \ \  1 \end{array}\right)$

\smallskip

- Le vecteur longueur associé est  : \ $\lambda = (\alpha,1-\alpha,
\alpha,1-\alpha, ....,\alpha,1-\alpha)$.

\smallskip

- Les discontinuités de \ $T$ \ sont : $\alpha, 1, 1+\alpha, 2 ....,
(n-1)+\alpha$.

\smallskip

- Les translations de \ $T$ \ sont : \ $\delta_i =  (n-i +1) -
\alpha$, \ $1\leq i\leq m$.

\smallskip

- Les images des discontinuités de \ $T$ \ sont : $0$, $1-\alpha,
1 , 2-\alpha, 2,....,n-\alpha$.

\bigskip

\bigskip

{\unitlength=0,28mm \hskip 2 truecm \begin{picture}(300,300)

\put(0,0){\line(1,0){300}} \put(0,300){\line(1,0){300}}
\put(0,0){\line(0,1){300}} \put(300,0){\line(0,1){300}}

\put(20,0){\dashbox(0,300){}} \ \put(60,0){\dashbox(0,300){}}
\put(80,0){\dashbox(0,300){}} \ \put(120,0){\dashbox(0,300){}}

\put(260,0){\dashbox(0,300){}} \ \put(240,0){\dashbox(0,300){}}
\put(200,0){\dashbox(0,300){}} \

\put(0,40){\dashbox(300,0){}}
 \put(0,60){\dashbox(300,0){}}
\put(0,100){\dashbox(300,0){}}
 \put(0,120){\dashbox(300,0){}}

\put(0,280){\dashbox(300,0){}}
 \put(0,240){\dashbox(300,0){}}
\put(0,220){\dashbox(300,0){}}

\put(-10,-10){$\scriptstyle 0$} \put(18,-10){$\scriptstyle {\alpha}$}
  \put(57,-10){$\scriptstyle 1 $}
\put(74,-10){$\scriptstyle 1 + \alpha$}
 \put(118,-10){$\scriptstyle 2$}
\put(227,-12){${ \scriptstyle n-1} $}
\put(256,-12){$\scriptstyle n-1+\alpha$}
\put(305,-12){$\scriptstyle n$}

 \put(-35,40){$\scriptstyle 1-\alpha$}
 \put(-15,60){$\scriptstyle 1$}
 \put(-35,100){$\scriptstyle 2-\alpha$}
 \put(-15,120){$\scriptstyle 2$}

\put(-15,300){$\scriptstyle  n$} \put(-35,280){$ \scriptstyle  n-\alpha$}
\put(-35,240){$\scriptstyle  n-1$}
\put(0,280){\line(1,1){20}}

\put(20,240){\line(1,1){40}}

\put(60,220){\line(1,1){20}}

\put(200,60){\line(1,1){40}}

\put(240,40){\line(1,1){20}}

\put(260,0){\line(1,1){40}}

 \end{picture}}\\\\
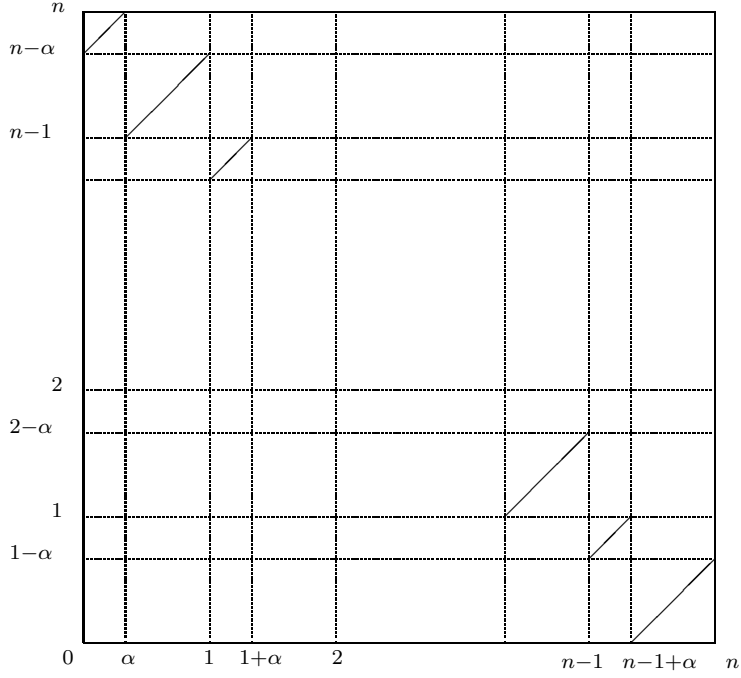
\centerline{\it Figure 5}

\bigskip

\bigskip

Pour \ $n\geq 3$, \ $T$ \ n'est pas minimal car \
$T([0,1]\cup[n-1,n]) = [0,1]\cup [n-1,n]$, aussi car \ $\sigma$ \
n'est pas un \ $n$-cycle. D'autre part, puisque \ $\delta_i = -
\alpha + (n-i +1)$, $1\leq i\leq m$, \ donc $T$ vérifie les
conditions de la proposition 2.1.A) avec \ $r = -\alpha $ et $s =
1$. Par suite, $f(x) = \exp(2i\pi x)$ est une fonction propre
associ\'ee \`a la valeur propre $\exp(- 2i\pi \alpha)$. 

\bigskip

{\bf Remarque.} Dans cette famille $(T_\alpha)$  à un paramètre \ $\alpha$ \
d'échanges d'intervalles de rang \ $2$, \ on voit que le non mélange
faible topologique et l'existence d'une valeur propre rationnelle non triviale
$exp (2i\pi \frac {1}{n})$ sont des  propriétés stables, comme conjecturé
dans  \cite{mdB88} .

\bigskip

$\bullet$ \ Cas  \ $\bf n = 2$: Echanges de 4 intervalles.

Dans ce cas, la permutation \ $\sigma = (21)$\ est un \ $2$-cycle,
donc par le lemme 4.1 \ $T$ \ est uniquement ergodique.
 Par contre,  \ $T$ \ ne vérifie pas la condition idoc car \ $T(\alpha) =
 \frac{1}{2}$. En divisant les longueurs par \ $2$ \ afin de se placer sur \ $I=
[0,1[$ \  et en changeant \ $\frac{\alpha}{2}$ \ par \ $\alpha$, \
on obtient :

\smallskip

 $\pi = \left (\begin{array}{ll}  1 \ \ 2 \ \ 3 \  \  4 \\  4 \ \ 3 \ \ 2 \  \
 1\end{array}\right)$ \ \ \ et \ \ \  $\lambda = ( \alpha,\frac{1}{2}-\alpha, \alpha,\frac{1}
 {2}-\alpha).$ D'o\` u
$$T(x)=\left\{
\begin{array}{llll}
x+1-\alpha & \textrm{si $x\in[0,\alpha[$}\\
x+\frac{1}{2}-\alpha &\textrm{si $x\in[\alpha,\frac{1}{2}[$}\\
x-\alpha&\textrm{si $x\in[\frac{1}{2},\frac{1}{2}+\alpha[$}\\
x-\alpha-\frac{1}{2}&\textrm{si $x\in[\frac{1}{2}+\alpha,1[$}
\end{array}
\right.$$

\medskip

\hskip 3 truecm
 \begin{picture}(180,180)

\put(0,0){\line(1,0){180}} \put(0,180){\line(1,0){180}}
\put(0,0){\line(0,1){180}} \put(180,0){\line(0,1){180}}

\put(-15,-15){$0$}

 \put(90,-20){${\frac{1}{2}}$}

 \put(20,-20){${\alpha}$}

 \put(110,-20){$\frac{1}{2}+ \alpha$}

  \put(180,-20){$1$}

 \put(-35,70){$\frac{1}{2}-\alpha$}

 \put(-15,90){$\frac{1}{2}$}

 \put(-35,160){$1-\alpha$}

 \put(-15,180){$1$}

\put(20,0){\dashbox(0,180){}} \put(90,0){\dashbox(0,180){}}
\put(110,0){\dashbox(0,180){}}

\put(0,70){\dashbox(180,0){}} \put(0,90){\dashbox(180,0){}}
\put(0,160){\dashbox(180,0){}}

\put(0,160){\line(1,1){20}}

\put(20,90){\line(1,1){70}}

\put(90,70){\line(1,1){20}}

\put(110,0){\line(1,1){70}}

 \end{picture}\\\\

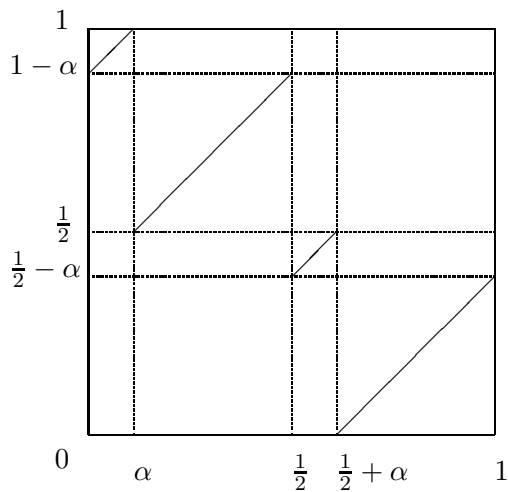
\centerline{\it Figure 6}

\bigskip

On vérifie directement que \ $f(x) = \exp(4i\pi x)$ \ est une
fonction propre de \ $U_{T}$ \ de valeur propre \ $\exp(-4i\pi
\alpha)$ et que $\exp(i\pi )$  est valeur propre
puisque $T^2([0, \frac{1}{2}[) =[0, \frac{1}{2}[$.

\bigskip

\bibliographystyle{amsplain}

\vfill

\textbf{Remerciements}. Je tiens \`a exprimer ma reconnaissance
aux professeurs Isabelle Liousse et Habib Marzougui pour leurs
encouragements et leurs aides.

\vfill
\end{document}